\documentclass[12pt,article]{amsart}

\usepackage{amsmath,amssymb}
\usepackage{tikz} 
\usepackage[margin=1in]{geometry}
\usepackage{amsmath,amssymb}
\usepackage{tikz}
\usepackage[normalem]{ulem}
\usepackage{hyperref}
\usepackage{amsthm}
\usepackage{cleveref}
\newtheorem{theorem}{Theorem}
\newtheorem{lemma}{Lemma}

\theoremstyle{definition}
\newtheorem{definition}[theorem]{Definition}
\newcommand{\Cay}{\operatorname{Cay}}
\newcommand{\low}{\operatorname{low}}
\newcommand{\Cp}{C'(1/8)}  
\newcommand{\abs}[1]{\left\lvert #1 \right\rvert}
\newcommand{\Rel}{\mathcal{R}} 

\title[Diagram of Geodesic Intervals in C'(1/8)--Groups]{Diagram of Geodesic Intervals in C'(1/8)--Groups}
\author{Alexander Hoffmann}

\begin{document}
\begin{abstract}
    Let $\Gamma$ be a group characterized by a small cancellation $C'(1/8)$-presentation. We study the structure of its Cayley graph, specifically geodesic intervals and digons. We prove that all geodesic segments [x,y] of length $n$ can be covered by finitely many simple geodesic digons with pairwise $C'(1/8)$-short lower intersections whose lower sides lie on a common geodesic segment [x,y]. This is a generalization of a proposition of Kharlampovich-Sklinos for random groups. We then bound this resulting diagram of digons, showing that the total perimeter of its cells and its edge count is bounded by $n$, and the cycle rank is equal to the number of distinct cell contours. Thus, although there may be many different geodesics, the size of their union is bounded linearly in the distance between their endpoints. The constant is independent of the number and lengths of the defining relators. In particular, the estimate applies to infinite presentations with unbounded relator lengths. From this we get a consequence for growth, allowing us to compare the number of geodesic words with the number of group elements of the same length.
\end{abstract}
\maketitle

\section{Introduction}\label{sec:intro}
Small cancellation theory studies groups whose defining relators have little overlap with each other. This condition strongly dictates the group's Van Kampen diagram and the geometry of the Cayley Graph \cite{LS}. A universal question about any finitely generated group is how its geodesics, the shortest path between two points, behave: the number of distinct geodesics between two given points, how far apart these geodesics can be, and how the set of all geodesics from a point $x$ to a point $y$ sit inside the Cayley Graph. For groups with a $\Cp$-presentation, the geometry of simple geodesic digons and simple geodesic triangles has been classified by Strebel \cite{Str}

In this paper we study the union of all geodesics between two points of the Cayley graph of a $\Cp$-group. Our starting point is a proposition of Kharlampovich and Sklinos for $\Cp$-groups with all relators being the same length \cite{KS}, which desribes geodesic intervals by digons with short lower intersections. We generalize this statement to all $\Cp$-groups where the relators can differ in lenght. Specifically, in \cref{thm:cover} we show that for any geodesic segment $[x,y]$ there are finitely many digons $D_1,\dots,D_k$ whose
lower sides lie on $[x,y]$ and pairwise have $\Cp$-short intersections, such that every geodesic from $x$ to $y$ is contained in $[x,y]\cup D_1\cup\dots\cup D_k$. The key step is a merging argument (\cref{lem:union}): two digons whose lower sides intersect in a long piece must share cells, and their union forms another digon. So the digons can be reduced until all pairwise intersections are short, at which point \cref{lem:short} shows that no further merging is possible. 

Our second result, \cref{thm:bound}, shows that the resulting diagram is small. Writing $n=d(x,y)$, we prove that the total perimeter of the distinct cells is less than $\tfrac83 n$, that the subgraph $H$ formed by $[x,y]$ and the digons has fewer than $\tfrac73 n$ edges, and that its cycle rank equals the number $N$ of distinct cell contours. In other words, although there could be many geodesics from $x$ to $y$, their union is bounded linearly by $n$. The proof combines the lower bounds of arcs from \cref{lem:arcs} along with a counting argument for the divisors and overlaps. A height function $h(z)=d(x,z)$ shows that nonconsecutive cells of the diagram are disjoint in the Cayley graph. We also prove that the number of geodescis is bounded by $N$, that $|\mathcal G(x,y)|\le 2^N$, by creating a tree and showing that a geodesic is determined uniqeuly by the edges outside of the tree that it uses. The constants do not depend on the number of relators or their lengths, so the bounds hold for infinite presentations with unbounded relator lengths. If every relator has length at least $m$, the bound for the number of distinct cell contours is improved, to $N<8n/(3m)$

Finally, we use compare two different growth rates. Let $a_n$ be the number of
group elements of length $n$ and $b_n$ the number of geodesic words of length $n$, where distinct words representing the same element are counted separately. We use \cref{thm:bound} to show $a_n\le b_n<2^{8n/(3m)}a_n$, and therefore $0\le h_{\mathrm{geo}}-h_\Gamma\le 8\log 2/(3m)$ (\cref{thm:growth}). Thus in a $\Cp$-group whose relators have a minimum length, the exponential grwoth rate of geodesic words is close to the group itself. 

\medskip
Geodesics in infinitely presented $\Cp$-groups were also studied by Arzhantseva and
Dru\c{t}u \cite{AD}. One of their technical results \cite[Theorem~3.15]{AD}
describes, for any two vertices $a,b$ of the Cayley graph, a set containing all
geodesics from $a$ to $b$. This set is a finite alternation of
\emph{compulsory} geodesic segments, on which the geodesic between the endpoints
is unique. They further use this to construct quasi-homomorphisms and to study bounded cohomology. 
Our \cref{thm:cover} is in the same spirit as \cite[Theorem~3.15]{AD}: both
place all geodesics between two points inside an explicit union of cells.
In Arzhantseva and Dru \cite{AD}, the set is organized along the geodesic, as
compulsory segments 
alternating with successions of contours, and the emphasis is on the
decomposition itself and its algebraic form (every group element is a unique
product of compulsory and relator-covered elements). We instead organize the set
as a finite family of digons based on a fixed geodesic $[x,y]$ with pairwise
$\Cp$-short lower intersections, following the random-group statement of
\cite{KS}, and we then quantify it in \cref{thm:bound} and \cref{thm:growth}. Conversly, we do not address the algebraic decomposition of group elements or the applications to quasi-homomorphisms
developed there.

\section{Preliminaries}

Throughout, we work with the group $\Gamma=\langle\Sigma_\Gamma\mid\Rel\rangle$, where 
$\Sigma_\Gamma$ is a finite set of generators and $\Rel$ is a (possibly infinite) set 
of relators. We also write $S=\Sigma_\Gamma$ and work with words over
$S^{\pm1}=S\cup S^{-1}$. The set $\Rel$ is \emph{symmetrized} if every $r\in\Rel$ is cyclically  reduced and $\Rel$ is closed under taking inverses and cyclic permutations. This group has the $\Cp$ condition which works as follows.

\begin{definition}[$\Cp$]\label{def:cp}
A word $p$ is a \emph{piece} (with respect to $\Rel$) if there are distinct
$r_1,r_2\in\Rel$ with $r_1=pu_1$ and $r_2=pu_2$ as words. The group satisfies \emph{$\Cp$} if every piece $p$ that is an initial subword of some $r\in\Rel$ satisfies $\abs{p}<\tfrac18\abs{r}$. In the Van Kampen diagram, we will think of this condition as that any two cells can not overlap by more than $\frac18$ of each cell's lengths. 
\end{definition}

\begin{definition}[Cayley Graph and Van Kampen diagram]\label{def:Cay and VK}
    The \emph{Cayley graph} $X=\Cay(\Gamma,\Sigma_\Gamma)$ has vertex set $\Gamma$ and
an edge of length $1$ joining $g$ to $gs$ for each $g\in\Gamma$ and
$s\in S^{\pm1}$. This induces a path metric $d$ on $X$ and $\abs{g}_S=d(1,g)$ for the word length of $g$. A \emph{geodesic segment} $[x,y]$ is an edge path from $x$ to $y$ of length $d(x,y)$, or a path with the shortest possible length. A \emph{Van Kampen diagram} over $\Gamma$ is a finite, connected, simply connected planar $2$-complex $D$ whose edges are labelled by $S^{\pm1}$ and such that the label of the boundary of each face, read from any vertex in either direction, is a relator in $\Rel$. The faces are called \emph{cells}. For a cell $C$ we write $\partial C$ for its boundary path or contour and $\abs{C}$ for its length, which is the length of the relator labelling it.
\end{definition}

\begin{definition}[Simple geodesic digon]\label{def:digon}
    Let $\Gamma = \langle \Sigma_\Gamma \mid  R\ \rangle$ be given by a symmetrized
$C'(1/8)$-presentation. A \emph{simple geodesic digon} with vertices $x,y$ in a Van Kampen Diagram $D$ is the union of two distinct geodesic between $x$ and $y$ such that the two geodesics only intersect at these points. The two sides are labelled as $\low(D)$ and $\operatorname{up}(D)$. The designations of the lower and upper sides of the digon are arbitrary and serve only as labels; they do reflect any intrinsic distinction in the digon's presentation. Strebel has classified the shape of simple geodesic digons in $\Cp$-groups \cite{Str}; they are divided into a chain of cells split by edges going from $\low(D)$ to $\operatorname{up}(D)$. We call the points where the divisor intersects the sides of the digon as division points. 
\end{definition}

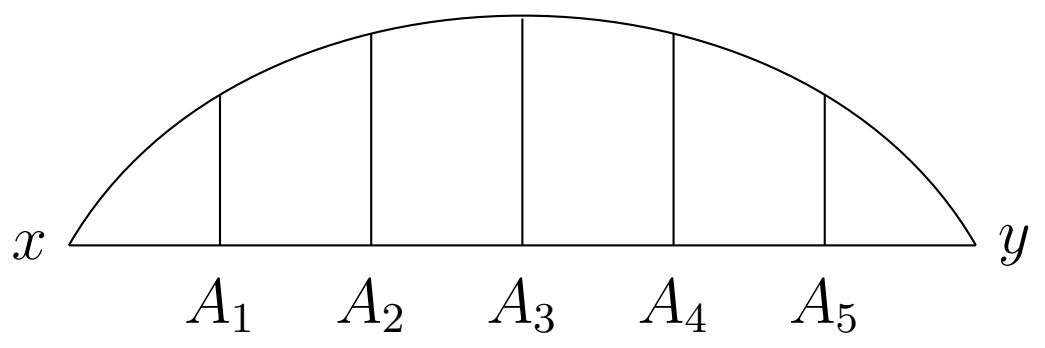
\begin{figure}[ht]
\centering
\begin{tikzpicture}[scale=1.0]
  \draw (0,0) -- (6,0) node[right] {$y$};
  \node[left] at (0,0) {$x$};
  \draw (0,0) to[out=60,in=120] (6,0);
  \draw (1,0) -- (1,1);
  \draw (2,0) -- (2,1.4);
  \draw (3,0) -- (3,1.5);
  \draw (4,0) -- (4,1.4);
  \draw (5,0) -- (5,1);
  \node at (1,-0.4) {$A_1$};
  \node at (2,-0.4) {$A_2$};
  \node at (3,-0.4) {$A_3$};
  \node at (4,-0.4) {$A_4$};
  \node at (5,-0.4) {$A_5$};
\end{tikzpicture}
\caption{A simple geodesic digon $D$ between $[x,y]$ with lower side $\low(D)$ and upper side $\operatorname{up}(D)$. Filled with cells and having lower division points $A_i$}
\label{fig:digon}
\end{figure}

For the remainder of this paper, the term \emph{digon} will be used as a shorthand for \emph{simple geodesic digon} unless stated otherwise.

\begin{definition}[$\Cp$-short intersection]\label{def:short}
For a cell $C$ of $D$, let
$\low_D(C)$ denote the maximal subpath of $\partial C$ contained in $\low(D)$.
Let $D,E$ be digons whose lower sides lie on the same geodesic segment. Their
lower sides have \emph{$\Cp$-short intersection} if for every cell $C$ of $D$,
every cell $C'$ of $E$, and every common subpath
$P \subseteq \low_D(C)\cap\low_E(C')$ we have
\[
  \abs{P} < \tfrac{1}{8}\min\{\abs{C},\abs{C'}\}.
\]
\end{definition}
\begin{figure}[ht]
\centering
\begin{tikzpicture}[scale=0.64]
  \draw (0,0) -- (6,0);
  \node[left] at (0,0) {$x$};
  \node[right] at (11.5,0) {$y$};
  \draw (0,0) to[out=60,in=120] (6,0);
  \draw (1,0) -- (1,1);
  \draw (2,0) -- (2,1.4);
  \draw (3,0) -- (3,1.5);
  \draw (4,0) -- (4,1.4);
  \draw (5,0) -- (5,1);
\begin{scope}[shift={(5.5,0)}]
 \begin{scope}[yscale=-1]
  \draw (0,0) -- (6,0);
  \draw (0,0) to [out=60,in=120] (6,0);
  \draw (1,0) -- (1,1);
  \draw (2,0) -- (2,1.4);
  \draw (3,0) -- (3,1.5);
  \draw (4,0) -- (4,1.4);
  \draw (5,0) -- (5,1);
  \end{scope}
  \end{scope}
\end{tikzpicture}
\caption{Two geodesic digons with lower side on geodesic $[x,y]$ having $\Cp$-short intersection}
\label{fig:two digons}
\end{figure}
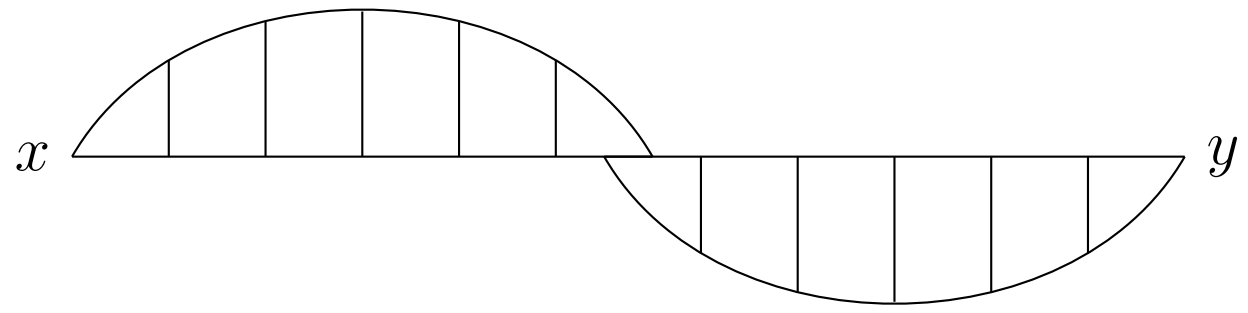

\section{Digons along common geodesic}

\begin{lemma}\cite[Lemma~6.9]{KS}\label{lem:arcs}
Let
\[
R=[u,c,v,d]
\]
be the boundary of a cell $C$ which is part of a digon
in $\Gamma$, with sides c and d being shared with other cells of the digon. Then the lengths of $u$ and $v$ are strictly greater than $\frac14 |C|$.
\end{lemma}

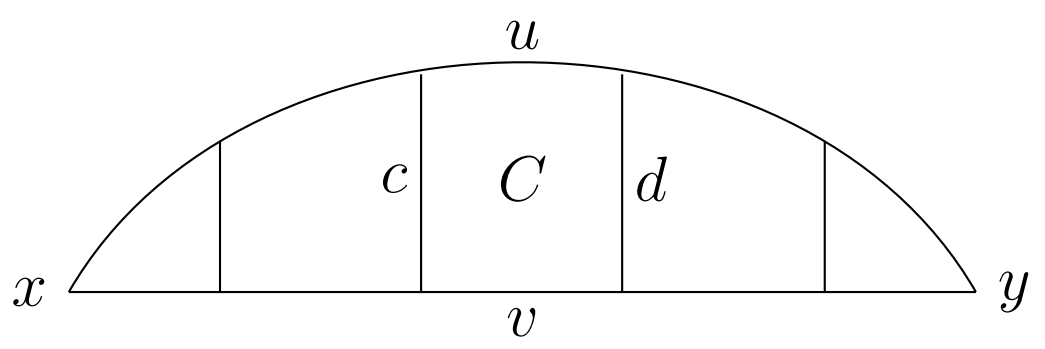
\begin{figure}[ht]
\centering
\begin{tikzpicture}[scale=1.0]
  \draw (0,0) -- (6,0) node[right] {$y$};
  \node[left] at (0,0) {$x$};
  \draw (0,0) to[out=60,in=120] (6,0);
  \draw (1,0) -- (1,1);
  \draw (2.33,0) -- (2.33,1.44);
  \draw (3.66,0) -- (3.66,1.44);
  \draw (5,0) -- (5,1);
  \node at (3,0.75) {$C$};
  \node at (3,1.7) {$u$};
  \node at (3,-0.2) {$v$};
  \node at (2.15,0.75) {$c$};
  \node at (3.85,0.75) {$d$};
\end{tikzpicture}
\caption{A geodesic digon $D$ between $[x,y]$ filled with cells having lower side $\low(D)$ and upper side $\operatorname{up}(D)$.}
\label{fig: lemma 1 digon}
\end{figure}
\begin{proof}
Because $c,d$ are shared between cells, due to $C'(1/8)$-short
cancellation, $|c|,|d|<\frac18|C|$. Using this we get
\begin{equation}\label{eq:lemma1}
    |u|+|v|=|C|-|c|-|d| > \frac34|C|
\end{equation}

\vspace{1em}

Because $u$ is a geodesic segment, the shortest
path between two points, $|u|,|v| \leq \frac12|C|$, and from \eqref{eq:lemma1}, we get $|u|,|v|>\frac14|C|$
\end{proof}

\begin{lemma}\cite[Lemma~6.15]{KS}\label{lem:union}
Suppose $D_1,D_2$ are digons for which $\operatorname{low}(D_1)\cup\operatorname{low}(D_2)=[x,y]$ is a geodesic segment and $D_1,D_2$ does not have
$C'(1/8)$-short intersection. Then for some digons $D$, $E$
\[
\operatorname{low}(D_1)\cup\operatorname{low}(D_2)
=
\operatorname{low}(D) \;\text{and} \;
\operatorname{low}(D_1)\cap\operatorname{low}(D_2)
=
\operatorname{low}(E).
\]
\centering And from this follows that $D_1\cup D_2 \subseteq D$
\end{lemma}
\begin{proof}
Let $A_1,\ldots,A_n$ be the division points of $\operatorname{low}(D_1)$ from the cells of $D_1$, and now add $A_0$ and $A_{n+1}$ to be the initial and terminal points of
$\operatorname{low}(D_1)$, with digon $D_1$ having $n+1$ cells. 
Likewise, let $A_0',A_1',\ldots,A_{m+1}'$ be the division, initial and terminal points of $D_2$, and note that 
\[
x=A_0,
\qquad
y=A_{m+1}'.
\]

\vspace{1em}

We first claim that the first cell of $D_2$ coincides with some
cell of $D_1$.

If $A_0'\in (A_n,A_{n+1}),$ then the first cell of $D_2$ and the last cell of $D_1$
would not have $C'(1/8)$-short intersection, so the piece of each
cell on $[x,y]$ would be greater than $\frac18\min\{|C|,|C'|\}$, so
\[
[A_0',A_{n+1}]
>
\frac18\min\{|C|,|C'|\},
\]
where $C$ is the last cell of $D_1$ and $C'$ is the first
cell of $D_2$.

However, because we are in a $C'(1/8)$ small cancellation group,
no two cells can overlap by more than
$\frac18\min\{|C|,|C'|\}$, so these two cells must coincide.

\vspace{1em}

If $A_0'\in[A_k,A_{k+1}]$, where $k<n$, then from \cref{lem:arcs},
\[
[A_j,A_{j+1}]>\frac14|C| \qquad [A_0',A_1']>\frac14|C'|
\]
for all $0\leq j\leq n$, where $C$ is the cell containing $[A_j,A_{j+1}]$ and $C'$ is the cell containing $[A'_0,A'_1]$. From this we get
\[
[A_0',A_1']
\cap
[A_j,A_{j+1}]
>
\frac18\min\{|C|,|C'|\}
\]
for either $j=k$ or $j=k+1.$
Consequently, the first cell of $D_2$ corresponds with the
$(j+1)$th cell of $D_1$.

\vspace{1em}

Let the first cell of $D_2$ coincide with the $(j+1)$th cell
of $D_1$.

We prove that the rest of the cells of $D_2$ coincide with cells of $D_1$ by induction on $r$, where $r$ denotes the cell number in $D_2$.
The case $r=1$ was established above: the first cell of $D_2$ coincides with the $(j+1)$th cell of $D_1$.

Now suppose that, for some $r \geq 1$, the $r$th cell of $D_2$ coincides with the $(j+r)$th cell of $D_1$, and suppose that both the $(r+1)$th cell of $D_2$ and the $(j+r+1)$th cell of $D_1$ exist. Since $\operatorname{low}(D_1) \cap \operatorname{low}(D_2)$ is a geodesic segment and $A'_r = A_{j+r}$, we get
\[
[A'_r,A'_{r+1}] \cap [A_{j+r},A_{j+r+1}]
=
[A'_r,A'_{r+1}]
\quad\text{or}\quad
[A_{j+r},A_{j+r+1}],
\]
since one of these segments contains the other.

By \cref{lem:arcs}, both segments have length at least $\frac{1}{4}|C|$, where $C$ is the corresponding cell. Therefore, their intersection has length at least $\frac{1}{8}|C|$.
By the $C'(1/8)$ small-cancellation condition, the $(r+1)$th cell of $D_2$ coincides with the $(j+r+1)$th cell of $D_1$.

Thus, by induction, the $r$th cell of $D_2$ coincides with the $(j+r)$th cell of $D_1$ for every $r$ for which these cells exist.

\vspace{1em}

Now we show the union of the digons forms a digon.

Suppose $j\neq n$, so the overlap is not just the final cell of $D_1$.

Let $B_i$ be the division points on $\operatorname{up}(D_1)$, and
similarly let $B_i'$ be the division points on $\operatorname{up}(D_2)$,
with the division points labelled from $0$ to $n+1$ and $0$ to $m+1$ respectively.

Note that
\[
B_{j+1}\in\operatorname{up}(D_2)
\]
(when $j=n-1$, $B_{j+1}$ should be considered to be $A_{n+1}$),
otherwise
\[
B_{j+1}\in[A_1',B_1'] \quad\text{or}\quad B_{j+1}\in\operatorname{low}(D_2).
\]

In both cases,
\[
[B_j,B_{j+1}]
\cup
[B_{j+1},B_{j+2}]
\]
is not a geodesic, a contradiction.

Therefore,
\[
B_{j+1}\in\operatorname{up}(D_2) \quad\text{and}\quad B_{j+1}=B_1'.
\]

Now we show the union
\[
[A_0,B_{j+1}]\cup[B_{j+1},A_{m+1}'],
\]
where
\[
(A_0,B_{j+1})\subseteq \operatorname{up}(D_1) \quad\text{and}\quad (B_{j+1},A_{m+1}')\subseteq \operatorname{up}(D_2)
\]
is a geodesic. With this the two geodesic segments
\[
\operatorname{low}(D_1)\cup\operatorname{low}(D_2) \quad\text{and}\quad (A_0,B_{j+1})\cup(B_{j+1},A_{m+1}')
\]
define a digon.

Note that
\[
d(x,y)
=
|[A_0,A_0']\cup[A_0',A_{m+1}']|
=
|[A_0,A_0']\cup[A_0',B_{j+1}]
\cup[B_{j+1},A_{m+1}']|,
\]
since
\[
[A_0',B_{j+1}]\cup[B_{j+1},A_{m+1}']
=
\operatorname{up}(D_2) \; \text{is a geodesic}
\]
Also
\[
\left|[A_0,A_0']\cup[A_0',B_{j+1}]\right|
\geq
\left|[A_0,B_{j+1}]\right|,
\]
so
\[
[A_0,B_{j+1}]\cup[B_{j+1},A_{m+1}'] \; \text{is a geodesic}
\]

\vspace{1em}

Suppose $j=n$. Now that the overlap begins in the last cell of $D_1$,
the point $B_{j+1}$ doesn't exist as an interior subdivision point.

The point $B_1'$ must be in $(A_{n+1},B_n)$, strictly less than $\frac18|C|$ away from $A_{n+1}$, where $C$ is the smaller cell between the last $D_1$ cell and the first $D_2$ cell.

We need to show that $[A_0,B_1']\cup[B_1',A_{m+1}']$ is a geodesic segment.
From digon $D_2$ we get
\[
|[A_0',B_1']|
+
|[B_1',A_{m+1}']|
=
|[A_0',A_{n+1}]|
+
|[A_{n+1},A_{m+1}']|,
\]
so
\[
|[A_0,A_{m+1}]|
=
|[A_0,A_0']|
+
|[A_0',B_1']|
+
|[B_1',A_{m+1}']|.
\]

Also,
\[
|[A_0,B_1']|
\leq
|[A_0,A_0']|
+
|[A_0',B_1']|
\; \text{because} \; [A_0,B_1'] \; \text{is a geodesic}.
\]

Therefore,
\[
|[A_0,A_{m+1}]|
\geq
|[A_0,B_1']\cup[B_1',A_{m+1}']|,
\]
and
\[
[A_0,B_1']\cup[B_1',A_{m+1}']
\]
is a geodesic segment. Hence,
\[
\operatorname{low}(D_1)\cup\operatorname{low}(D_2)
\]
is a digon.

\vspace{1em}

Now we show the intersection of the digons forms a digon.

For $j=n$, $\operatorname{low}(D_1)\cap\operatorname{low}(D_2)
=
\operatorname{low}(E)$
for the digon $E$ with a single cell determined by the geodesic segments
\[
[A_0',A_{n+1}] \quad\text{and}\quad [A_0',B_n]\cup[B_n,A_{n+1}].
\]

We need to show that
\[
|[A_0',A_{n+1}]|
=
|[A_0',B_1']|
+
|[B_1',A_{n+1}]|.
\]

From digon $D_1$ we get
\begin{equation}\label{eq:D1 geodesic}
    |[A_0,A_{n+1}]|
=
|[A_0,B_n]|
+
|[B_n,A_{n+1}]|.
\end{equation}

Since $\operatorname{low}(D_1)\cup\operatorname{low}(D_2)$ is a side of a digon $D$, and $[B_n,A_{n+1}]$ is a geodesic segment containing $B_1'$, we get
\begin{equation}\label{eq:big geodesic}
    |[A_0,A_{n+1}]|
+
|[A_{n+1},A_{m+1}']|
=
|[A_0,B_n]|
+
|[B_n,A_{n+1}]|
-
|[A_{n+1},B_1']|
+
|[B_1',A_{m+1}']|.
\end{equation}

Combining \eqref{eq:D1 geodesic} and \eqref{eq:big geodesic}, we get
\begin{equation}\label{eq:combined geodesic}
    |[A_{n+1},A_{m+1}']|
=
|[B_1',A_{m+1}']|
-
|[A_{n+1},B_1']|.
\end{equation}

From digon $D_2$ we get
\begin{equation}\label{eq:D2 geodesic}
    |[A_0',B_1']|
+
|[B_1',A_{m+1}']|
=
|[A_0',A_{n+1}]|
+
|[A_{n+1},A_{m+1}']|,
\end{equation}
and combining \eqref{eq:combined geodesic} and \eqref{eq:D2 geodesic} we finally show
\[
|[A_0',A_{n+1}]|
=
|[A_0',B_1']|
+
|[B_1',A_{n+1}]|.
\]

The case for $j<n$ follows by a similar argument.
\end{proof}

\begin{lemma}\cite[Lemma~6.19]{KS}\label{lem:short}
Let $D_1,D_2$ be digons for which
\[
\operatorname{low}(D_1)\cup\operatorname{low}(D_2)=[x,y]
\]
is a geodesic segment and let $D_1,D_2$ have
$C'(1/8)$-short intersection.
Then there is no digon or semidigon a side of which contains the union $\operatorname{low}(D_1)\cup\operatorname{low}(D_2)$.
\end{lemma}
\begin{proof}
Suppose a digon $D$ existed that contained
$\operatorname{low}(D_1)\cup\operatorname{low}(D_2)$ as a side.
By Lemma 2 both $D_1,D_2$ are contained in $D$. Let the cells
$C_1$ of $D_1$ and $C_2$ of $D_2$ be the cells that intersect. 
They must both coincide with a cell of $D$ or coincide with
consecutive cells of $D$. Both are impossible due to
$C_1$ and $C_2$ intersecting only on a short segment on $[x,y]$.
\end{proof}

\begin{theorem}\cite[Proposition~6.20]{KS}\label{thm:cover}
Suppose that $[x,y]$ is a geodesic segment in
$\operatorname{Cay}(\Gamma,\Sigma_\Gamma)$. Then there exist finitely many
digons $D_1,D_2,\ldots,D_k$ with the following properties:

\begin{itemize}
    \item For each $i\leq k$, $\operatorname{low}(D_i)\subseteq [x,y].$

    \item For any two $i,j\le k$,
    $\operatorname{low}(D_i)$ and $\operatorname{low}(D_j)$ have
    $C'(1/8)$-short intersection.

    \item If $J$ is a geodesic segment between $x$ and $y$, then it is
    contained in
    \[
    [x,y]\cup\{D_1,D_2,\ldots,D_k\}.
    \]
\end{itemize}
\end{theorem}
\begin{proof}
There can only be finitely many geodesic segments $[x,y]$:
In $\operatorname{Cay}(\Gamma,\Sigma_\Gamma)$ with a finite generating set,
$d(x,y)=m$ and each vertex has $k$ neighbors, so the total number of edge paths of length $m$ is $k^m$:
\[
\#([x,y])\le k^m<\infty.
\]

Let $J_0,J_1,\ldots,J_m$ be a list of them with $J_0=[x,y]$.
The union $J_0\cup J_i$ is equal to the union of $[x,y]$ with finitely many digons $D_1,\ldots,D_m$, whose lower sides are contained in $[x,y]$.

Now consider digons obtained by $J_1,J_2$.
If two digons do not have $C'(1/8)$-short intersection, then by
Lemma 2 their union is a lower side of a digon and we replace them with this digon.
\[
D_1\cup D_2\subseteq D,
\]
\[
[x,y]\cup D_1\cup D_2\subseteq [x,y]\cup D,
\]
so no geodesic segments are lost when replacing the two digons by their union. 

If two digons have $C'(1/8)$-short intersection, by Lemma 3 no other digon can contain the union of the lower sides.

Now assume we have a set of digons whose union with $[x,y]$
contains
\[
J_0,\ldots,J_i.
\]

Consider the digons obtained by $J_{i+1}$.
For each digon, we either merge it with some existing one, and if it has $C'(1/8)$-short intersection with all other digons, the union of the lower sides can't be part of some other digon.
By induction our proof is finished.
\end{proof}

\medskip

\section{The subgraph and its numerical invariants.}
Write $\mathcal G(x,y)$ for the set of
geodesic edge paths from a vertex $x$ to a vertex $y$.
Its elements are distinguished by their edge sequences.
The \emph{geodesic interval} between $x$ and $y$ is the subgraph
\[
 \mathcal I(x,y)=\bigcup_{I\in\mathcal G(x,y)}I.
\]
Given finitely many digons $D_1,\ldots,D_k$ based on geodesic $I$,
let $D_j^{(1)}$ denote the image in $X$ of all the edges of
the diagram $D_j$, including its divisors, and put
\[
 H=I\cup\bigcup_{j=1}^k D_j^{(1)}.
\]
Let $\mathcal C$ be the set of distinct contour images of all
cells in these diagrams, and define
\[
 N=|\mathcal C|,
 \qquad P=\sum_{C\in\mathcal C}|C|.
\]
For a finite connected graph $H$, its \emph{cycle rank} is
\[
 \beta_1(H)=|E(H)|-|V(H)|+1,
\]
where $E(H)$ and $V(H)$ are its edge and vertex sets.
For $k=0$ these conventions give $H=I$ and $N=P=0$.

\begin{theorem}\label{thm:bound}
Let $x,y$ be distinct vertices of $X$, let $n=d(x,y)$, and fix
a geodesic $I=[x,y]$.
There exist finitely many digons
$D_1,\ldots,D_k$ based on $I$, with pairwise short lower
intersections, such that
\[
 \mathcal I(x,y)\subseteq H,
 \qquad
 P<\frac83n,
 \qquad
 |E(H)|<\frac73n,
 \qquad
 \beta_1(H)=N.
\]
The empty family is allowed.
If $m$ is a positive integer such that $|r|\ge m$ for every
$r\in R$, then
\[
 N<\frac{8n}{3m},
 \qquad
 |\mathcal G(x,y)|\le 2^N<2^{\,8n/(3m)}.
\]
\end{theorem}
\begin{proof}
Choose the digons $D_1,\ldots,D_k$ from Theorem 4. This gives
\[
\mathcal I(x,y)\subseteq H,
\]
as all geodesics would follow along these digons.

Recall \cref{lem:arcs}, that if a cell has perimeter $L$, with lower and upper arcs
$u,v$, then
\[
|u|,|v|>\frac{L}{4}.
\]
We extend this specifically for the first or last cells of a digon. Since they have only one divisor, we can rewrite \eqref{eq:lemma1} in \cref{lem:arcs} to get
\[
|u|+|v|=|C|-|c| > \frac78|C|
\]
\[
|u|,|v|>\frac{3L}{8}.
\]

\vspace{1em}

We will first measure the lengths and show
\[
p\leq \frac{8}{3}n,
\qquad
|E(H)|<\frac{7}{3}n.
\]

Let $\ell_i$ be the length of either side of $D_i$. Let $s$ be the sum
of the lengths of all divisors. Let $o$ be the sum of the lengths of
all positive overlaps between lower sides. Let $\ell$ be the length
of the union of lower sides.

First, by considering $\ell_j$ as the lower sides of $D_j$, we see that the difference between the sum of all the lower sides and the union of the lower sides is where the lower sides overlap. So, 
\begin{equation}\label{eq:lower sides sum}
    \sum_j \ell_j = \ell + o
\end{equation}

Since all lower sides of digons are contained by $I = [x, y]$ with $d(x,y) = n$, it follows that 
\begin{equation}\label{eq:l less than h}
    \ell \le h
\end{equation}

Note that $P = \sum_{c \in C} |c|$, which is the sum of the lengths of every cell. All the lower and upper arcs of cells are contained in the lower and upper sides of digons, giving the total sum $2 \sum_j \ell_j$. Each divisor is used by exactly two cells, so using \eqref{eq:lower sides sum} gives us

\begin{equation}\label{eq:P sum}
P = 2 \sum_j \ell_j + 2s = 2\ell + 2o + 2s = 2\ell + 2(s+o)
\end{equation}

Now since every divisor or overlap belongs to two cells, and each arc is shorter than $1/8$ of its perimeter, we get 
\begin{equation}\label{eq:s+o bound}
    2(s+o) < P/4 
\end{equation}

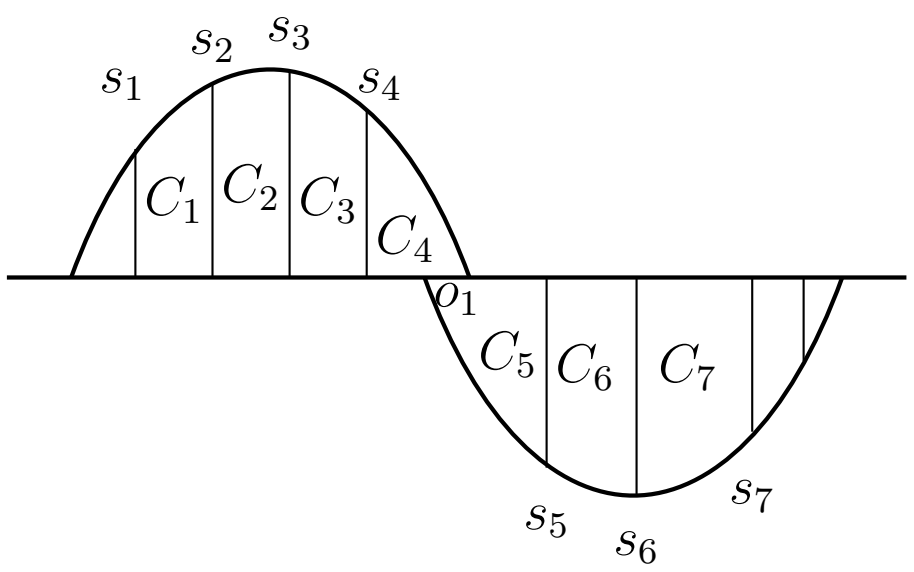
\vspace{0.5em}

\begin{figure}[ht]
\centering
\begin{tikzpicture}[scale=0.85, every node/.style={transform shape}]
  \draw[thick] (-0.5,0) -- (6.5,0);
  \draw[thick] (0,0) to[out=70,in=110,looseness=1.9] (3.1,0);
  \draw[thick] (2.75,0) to[out=290,in=250,looseness=1.9] (6,0);

  \draw (0.5, 0) -- (0.5, 1);
  \draw (1.1, 0) -- (1.1, 1.5);
  \draw (1.7, 0) -- (1.7, 1.6);
  \draw (2.3, 0) -- (2.3, 1.3);

  \node at (0.8, 0.6) {$C_1$};
  \node at (1.4, 0.7) {$C_2$};
  \node at (2, 0.6) {$C_3$};
  \node at (2.6, 0.3) {$C_4$};

  \node at (0.4, 1.5) {$s_1$};
  \node at (1.1, 1.8) {$s_2$};
  \node at (1.7, 1.9) {$s_3$};
  \node at (2.4, 1.5) {$s_4$};
  \node at (3.33, 0.1) [below left] {$o_1$};

  \draw (3.7, 0) -- (3.7, -1.48);
  \draw (4.4, 0) -- (4.4, -1.7);
  \draw (5.3, 0) -- (5.3, -1.2);
  \draw (5.7, 0) -- (5.7, -0.65);

  \node at (3.4, -0.6) {$C_5$};
  \node at (4.0, -0.7) {$C_6$};
  \node at (4.8, -0.7) {$C_7$};

  \node at (3.7, -1.9) {$s_5$};
  \node at (4.4, -2.1) {$s_6$};
  \node at (5.3, -1.7) {$s_7$};

  \node[anchor=west] at (7, 1.3) {$s_1 + s_2 < \frac{1}{4} |C_1|,\ s_2 + s_3 < \frac{1}{4} |C_2|$};
  \node[anchor=west] at (7, 0.5) {$s_3 + s_4 < \frac{1}{4} |C_3|,\ s_4 + o_1 < \frac{1}{4} |C_4|$};
  \node[anchor=west] at (7, -0.3) {$o_1 + s_5 < \frac{1}{4} |C_5|,\ s_5 + s_6 < \frac{1}{4} |C_6|$};

  \draw[->, thick] (9.5, -0.6) -- (9.5, -1.1);

  \node[anchor=west] at (6.8, -1.4) {$2(s_1 + \dots s_6 + o_1) < \frac{1}{4}(|C_1| + |C_2| \dots + |C_6|)$};
\end{tikzpicture}
\caption{A section showing how the divisors and overlaps are counted twice}
\label{fig:counting}
\end{figure}

\vspace{0.5em}

When inequalities for all cells are added up, the final equality becomes $2(s+o) < P/4$.

Combining \eqref{eq:l less than h}, \eqref{eq:P sum}, and \eqref{eq:s+o bound}, gives us a bound for $P$:
\begin{equation}\label{eq:P bound}
    P=2\ell+2(s+o),
\implies
P<2n+\frac{P}{4},
\implies
\frac{3P}{4}<2n,
\implies
P<\frac{8}{3}n.
\end{equation}

Now look at $|E(H)|$, the length of all edges. The length of the lower
edges is less than $n$, and from there we add the upper sides and
divisors. Then we use \eqref{eq:P sum} and \eqref{eq:P bound} to get
\[
|E(H)|\leq n+\sum_j\ell_j + s, 
\qquad
|E(H)|\leq n+\frac{P}{2}.
\]
\[
|E(H)|<\frac{7}{3}n.
\]

\vspace{1em}

We now show $\beta_1(H)=N,$

It follows from \cref{thm:cover} that each digon meets the interior of $I$ and the other
digons precisely in its lower side. Applying Euler's formula $V+F-E=2$
to a digon $D$ with $q$ cells, we get
\[
\beta_1(D)=|E(D)|-|V(D)|+1=q,
\]
since we don't include the outside of the digon as a face, unlike in Euler's formula.

Attaching this digon along its lower side increases the cycle rank of $H$
by $q$. So with $k$ digons with number of cells
$q_1,q_2,\ldots,q_k$,
\[
\beta_1(H)=q_1+q_2+\cdots+q_k.
\]
Since $N$ is the total number of cells,
\[
N=q_1+q_2+\cdots+q_k,
\]
and therefore
\[
\beta_1(H)=N.
\]

\vspace{1em}

Now we count the geodesics:

Choose a tree $T\subseteq H$ containing all vertices of $H$. $T$ is
connected and has no cycles. Notice that a tree containing $|V(H)|$
vertices has exactly $|V(H)|-1$ edges as each additional vertex after the first adds one edge. So the number of edges in $H$ but not in $T$ is
\[
|E(H)|-(|V(H)|-1) = N.
\]

We now show a simple path is determined by the set of edges outside of
$T$ that it uses. A simple path does not repeat vertices, so note that
\[
\#\{\text{geodesics}\}
\leq
\#\{\text{simple paths}\},
\]
since every geodesic is a simple path.

Suppose two distinct paths $P_1,P_2$ use exactly the same non-tree edges,
so the only place they differ is inside $T$. But if there were two
different paths in $T$ connecting the same vertices, it would give a
cycle, which $T$ cannot have. There is only one possible path through
$T$ between any two vertices, so $P_1$ and $P_2$ can't be distinct.

There are $N$ non-tree edges, giving us $2^N$ sets of non-tree edges. Therefore there are at most $2^N$ simple paths, and hence at most $2^N$ geodesic paths. 

Now if $m\leq |r|$ for every
$r\in R$ and some positive integer $m$, we get
\[
P=\sum_{i=1}^N |r_i|
=|r_1|+|r_2|+\cdots+|r_N|
\geq m+m+\cdots+m=mN,
\]
since there are $N$ cells and $P$ accounts for overlap between
cells. So
\[
P\geq mN
\qquad\text{and}\qquad
\frac{8n}{3m}\geq N,
\]
giving us the bound
\[
|\mathcal G(x,y)|\le 2^N<2^{\,8n/(3m)}.
\]
\end{proof}

\section{Growth Rates}

For a word $w$ over $S^{\pm1}$, let $\overline w$ denote the
element of $G$ it represents.
The word $w$ is \emph{geodesic} if $|w|=|\overline w|_S$.
For $n\ge1$, define
\[
 a_n=\bigl|\{g\in G:|g|_S=n\}\bigr|,
 \qquad
 b_n=\bigl|\{w\text{ over }S^{\pm1}:|w|=|\overline w|_S=n\}\bigr|.
\]
Thus $a_n$ counts group elements of length n, while $b_n$ counts geodesic words of length n;
distinct words representing the same element are counted
separately in $b_n$.
When $G$ is infinite, define the logarithmic exponential
growth rates, with respect to the fixed generating set, by
\[
 h_G=\limsup_{n\to\infty}\frac{\log a_n}{n},
 \qquad
 h_{\mathrm{geo}}=\limsup_{n\to\infty}\frac{\log b_n}{n},
\]
where $\log$ denotes the natural logarithm.

\begin{theorem}\label{thm:growth}
Suppose $G$ is infinite and $m$ is a positive integer such
that $|r|\ge m$ for every $r\in R$.
Then, for every $n\ge1$,
\[
 a_n\le b_n<2^{\,8n/(3m)}a_n.
\]
Consequently,
\[
 0\le h_{\mathrm{geo}}-h_G
       \le\frac{8\log 2}{3m}.
\]
\end{theorem}
\begin{proof}
Take $x=1$, $y=g$, and suppose $|g|_S=n$.
Then there are $a_n$ elements satisfying $|g|_S=n$,
each one having at most $2^{\,8n/(3m)}$ geodesic words from \cref{thm:bound}.

Every element $g$ has at least one geodesic word representing it,
so we get
\[
a_n\leq b_n\leq 2^{\,8n/(3m)}a_n.
\]

Now we take logarithms and divide by $n$ to get
\[
\frac{\log a_n}{n}
\leq
\frac{\log b_n}{n}
\leq
\frac{\log a_n}{n}
+\frac{8\log 2}{3m}.
\]

Now we take $\limsup$ and use our definitions of
$h_G$ and $h_{\mathrm{geo}}$ to get
\[
h_G\leq h_{\mathrm{geo}}
\leq h_G+\frac{8\log 2}{3m},
\implies
0\leq h_{\mathrm{geo}}-h_G
\leq \frac{8\log 2}{3m}.
\]
\end{proof}

\section*{Acknowledgements}

I would like to thank Professor O. Kharlampovich who introduced me to small cancellation groups and suggesting the problem that became the focus of the paper. Her guidance and encouragement were invaluable throughout this project.

\end{document}